\documentclass{article}
\usepackage{amssymb,amsmath}
\usepackage{graphics}

\def\qed{\hfill {\hbox{${\vcenter{\vbox{               
   \hrule height 0.4pt\hbox{\vrule width 0.4pt height 6pt
   \kern5pt\vrule width 0.4pt}\hrule height 0.4pt}}}$}}}

\def\tr{\triangleright}
\def\bar{\overline}

\newtheorem{definition}{Definition}

\newtheorem{example}{Example}
\newtheorem{remark}[example]{Remark}

{\qed\par\medskip}

\date{}

\title{\Large \textbf{Rack Module Enhancements of Counting Invariants}}

\author{\large
Aaron Haas${ }^{\dagger 1}$
\and 
Garret Heckel${ }^{\ast 2}$
\and 
Sam Nelson${ }^{\dagger 3}$
\and 
Jonah Yuen${ }^{\dagger 4}$
\and 
Qingcheng Zhang${ }^{\ddagger 5}$
} 

\begin{document}
\maketitle

\begin{abstract} We introduce a modified rack algebra $\mathbb{Z}[X]$ for 
racks $X$ with finite rack rank $N$. We use representations of 
$\mathbb{Z}[X]$ into rings, known as \textit{rack modules}, 
to define enhancements of the rack counting invariant for classical and 
virtual knots and links. We provide computations and examples to show that 
the new invariants are strictly stronger than the unenhanced counting 
invariant and are not determined by the Jones or Alexander polynomials.
\end{abstract}



\section{\large \textbf{Introduction}}

First introduced in \cite{AG}, a \textit{quandle module} is a representation
of an associative algebra $\mathbb{Z}[Q]$ known as the \textit{quandle 
algebra}, defined from a quandle $Q$. In \cite{CEGS} quandle modules were 
used to define knot and link invariants including some enhancements of the 
quandle counting invariant. In \cite{NJ} a generalized notion of rack modules 
defined in terms of trunks was introduced.

In this paper we introduce a modified form of the rack algebra for racks
with finite rack rank and use the resulting representations into modules
over finite rings to define new enhancements of the rack counting invariant 
$\Phi^{\mathbb{Z}}_{X}$ from \cite{N}. These new enhancements specialize to 
$\Phi^{\mathbb{Z}}_{X}$ but are able to distinguish many knots and links
with the same $\Phi^{\mathbb{Z}}_{X}$ values, and hence are strictly stronger
invariants. We provide examples which demonstrate that the new enhanced 
invariants can distinguish knots with the same Jones and Alexander polynomials.

The paper is organized as follows. In section \ref{rb} we review the basics
of racks and the rack counting invariant. In section \ref{ram} we define
the rack algebra and resulting rack modules. In section \ref{rinv} we use
finite rack modules to define new enhancements of the rack counting invariant.
In section \ref{comp} we give examples and describe our methods of 
computation of the new invariants. In section \ref{q} we collect a few
open questions and suggest directions for future work.

This paper is the end result of a summer 2010 project funded by Fletcher 
Jones Fellowships in Mathematics in the Claremont Center for Mathematical
Sciences, a summer undergraduate research program at the Claremont Colleges. 

\bigskip

\noindent
\parbox[b]{6.3in}{\rule{2cm}{0.2mm}
\footnotesize

${ }^{\dagger}$  \textsc{Department of Mathematics, Claremont McKenna College,
850 Columbia Ave, Claremont, CA 91711} 

${ }^{\ast}$  \textsc{Department of Mathematics, Claremont Graduate 
University, 150 E. 10th St., Claremont, CA 91711} 

${ }^{\ddagger}$ \textsc{Department of Mathematics, Pomona College, 
610 N. College Ave, Claremont, CA 91711} 

${ }^1$ahaas12@cmc.edu \
${ }^2$heckel.garret@gmail.com \ 
${ }^3$knots@esotericka.org \
${ }^4$jyuen11@cmc.edu\ 
${ }^5$qz002008@mymail.pomona.edu}

\newpage

\section{\large \textbf{Rack Basics}}\label{rb}

In this section, we review some rack basics needed for subsequent sections. 
We begin with a definition from \cite{FR}.
\begin{definition}
\textup{A \textit{rack} is a set $X$ with an operation 
$\tr: X \times X \rightarrow X$ satisfying the following axioms:
\begin{list}{}{}
\item[(i)] {for all $x,y \in X$, there exists a \textit{unique} 
$z \in X$ such that $x=z \tr y$, and} 
\item[(ii)] {for all $x,y,z \in X$, the equation 
\[(x \tr y) \tr z = (x \tr z) \tr (y \tr z)\] is satisfied.}
\end{list}
Note that Axiom (i) is equivalent to the following:
\begin{list}{}{}
\item[(i$'$)] {there exists a second operation 
$\tr^{-1}:X\times X \to X$ such that for all $x,y\in X$
\[(x \tr y) \tr^{-1} y = x=(x\tr^{-1} y)\tr y.\]}
\end{list}}
\end{definition}

As shown below, the operation $\tr$ denotes a crossing from right 
to left given the positive direction of the over-crossing strand. Conversely, 
the operation $\tr^{-1}$ denotes a crossing from the opposite direction.

\[\includegraphics{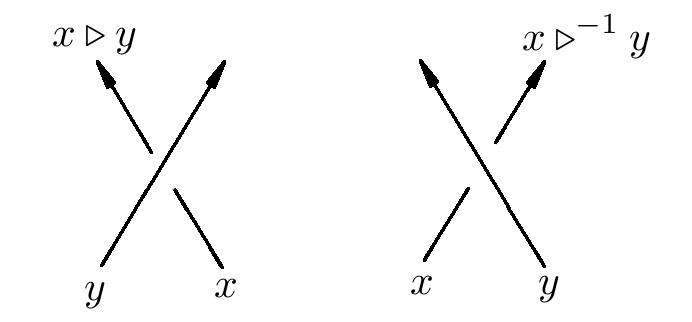}\]

If we regard $x, y \in X$ as two labels on arcs in a knot diagram, 
as shown below, we see that axioms (i) and (ii) correspond to 
Reidemeister moves II and III.

\[\scalebox{0.9}{\includegraphics{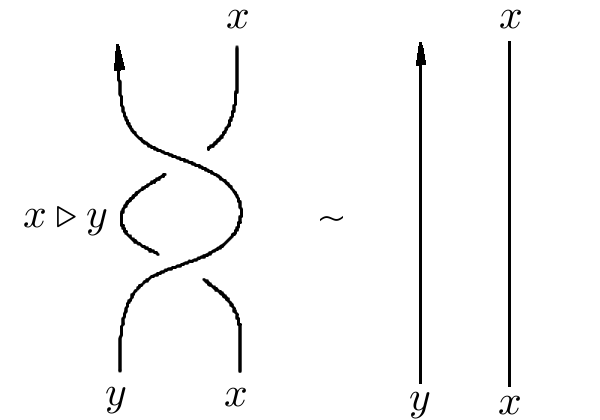}} \quad 
\scalebox{0.7}{\includegraphics{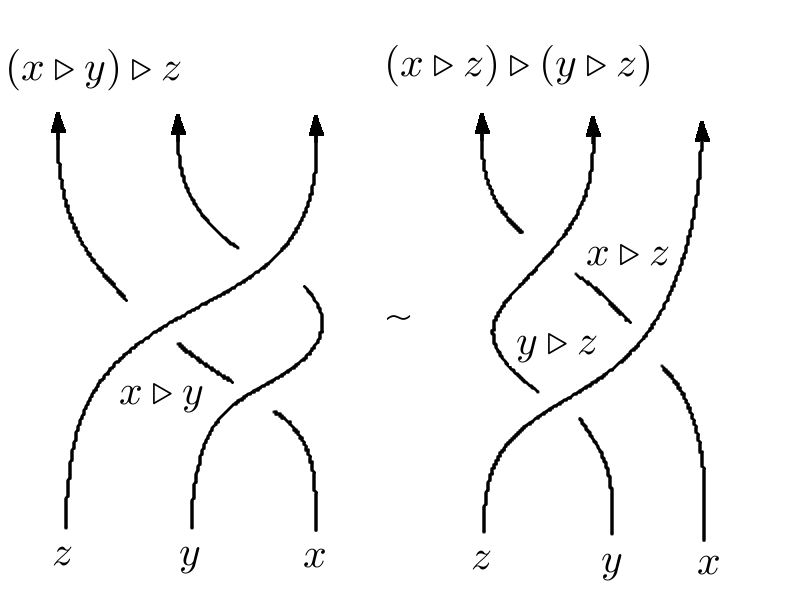}}\]

The equivalence relation on knot diagrams generated by the Reidemeister
II and III moves is known as \textit{regular isotopy}, and labelings of
links by racks are preserved by regular isotopy moves.
Rack colorings are also preserved by the \textit{blackboard framed type 
I moves}, which preserve the writhe or \textit{blackboard framing}:

\[\includegraphics{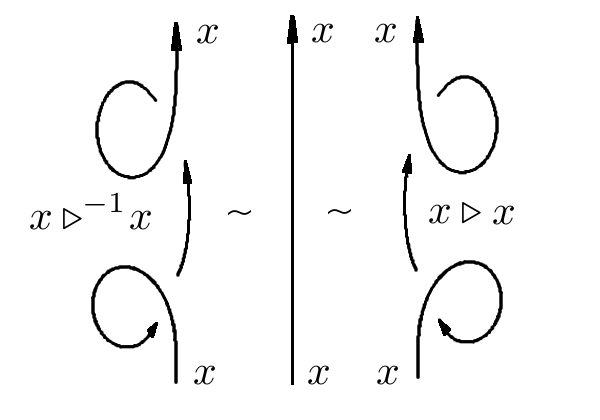}\]

\begin{example} \textup{\textit{Constant Action Racks}.
\\Given a set $X$, we can define a binary operation $x \tr y = \sigma(x)$, 
where $x,y \in X$ and the bijection $\sigma: X \rightarrow X$ is any 
permutation in the symmetric group $S_{X}$ on $X$. Then under the 
operation $\tr$, $X$ is a rack, since $\tr^{-1}=\sigma^{-1}(x)$ 
and} 
\[(x \tr y) \tr z = \sigma (x \tr y) = \sigma^2(x) 
= \sigma (x \tr z) = (x \tr z) \tr (y \tr z).\]
\end{example}

\begin{example} \textup{\textit{$(t,s)$-racks}.
\\A $(t,s)$-rack is a module $X$ over 
$\ddot{\Lambda} = \mathbb{Z}[t^{\pm1},s]/(s^2-(1-t)s)$ under the operations
\[ x \tr y  =  tx + sy \quad
\mathrm{and} \quad x \tr^{-1} y  = t^{-1}x - t^{-1}sy.\]
Convenient examples of $(t,s)$-rack structures include:
\begin{itemize}
\item $X=\mathbb{Z}_n$ with a choice of $s\in\mathbb{Z}_n$ and 
$t\in\mathbb{Z}_n^{\ast}$ satisfying $s^2=(1-t)s$, e.g. $\mathbb{Z}_4$ with
$t=1$ and $s=2$,
\item $X=\ddot\Lambda/I$ for an ideal $I$, e.g. $\ddot\Lambda/(t^2+1)$,
\item $X$ any abelian group with an automorphism $t:X\to X$ and an endomorphism
$s:X\to X$ satisfying $s^2=(Id-t)s=s(Id-t)$.
\end{itemize} 
If $s=1-t$ then our $(t,s)$-rack is a quandle, known as an \textit{Alexander 
quandle}.
}
\end{example}

\begin{definition}
\textup{Analogous to a Cayley table for a group, a \textit{rack matrix} 
describes an operation $\tr$ on a finite rack $X=\{x_1, x_2, \dots, x_n\}$. 
The entry $k$ of the rack matrix in the $i$th row and $j$th column is 
defined to be the subscript of $x_k$, where $x_k=x_i \tr x_j$.}
\end{definition}

\begin{example}
\textup{The following are the rack matrices of the constant action rack 
$X=\{x_1, x_2, x_3, x_4\}$ with $\sigma = (1234)$ and the 
$(t,s)$-rack $Y=\mathbb{Z}_4$ with $t=1$ and $s=2$:
\[M_X=
\left[\begin{array}{cccc}
2 & 2 & 2 & 2\\
3 & 3 & 3 & 3\\
4 & 4 & 4 & 4\\
1 & 1 & 1 & 1
\end{array}\right] \quad
M_Y=
\left[\begin{array}{cccc}
3 & 1 & 3 & 1 \\ 
4 & 2 & 4 & 2 \\ 
1 & 3 & 1 & 3 \\
2 & 4 & 2 & 4 
\end{array}\right].
\]}
\end{example}

Every tame blackboard-framed oriented knot or link $L$ has a 
\textit{fundamental rack} $FR(L)$ consisting of equivalence classes
of rack words in a set of generators corresponding one-to-one with
the set of arcs in a diagram of $L$ modulo the equivalence relation
generated by the rack axiom relations and the crossing relations in 
$L$. If we wish to specify explicitly the writhe vector or blackboard
framing vector of $L$, we will write $FR(L,\mathbf{w})$ where 
$\mathbf{w}=(w_1,\dots, w_n)$ and $w_k$ is the writhe of the $k$th
component of $L$.

For a given framed link $L$ and rack $X$, the number of rack labelings of
the arcs in $L$ by rack elements is an invariant of framed isotopy
known as the \textit{basic counting invariant}, denoted 
$|\mathrm{Hom}(FR(L),X)|$ since each labeling determines a unique rack 
homomorphism from $FR(L)$ to $X$.

\begin{definition}
\textup{Let $X$ be a rack. For every element $x \in X$, the \textit{rack rank} 
of $x$, $N(x)$, is the minimal integer $N$ such that $\pi^N(x)=x$
where the \textit{kink map} $\pi: X \rightarrow X$ is defined by 
$\pi(x)=x \tr x$. The \textit{rack rank} of $X$, $N(X)$, is the least 
common multiple (lcm) of the rank ranks of $x \in X$:}
\[N(X)=\mathrm{lcm}\{N(x) \ | \ x \in X \}.\] 
\textup{If $X$ is finite, then the rack rank is the exponent of the kink
map $\pi$ considered as an element of the symmetric group $S_X$.}
\end{definition}

Rack labelings of knot and link diagrams by a rack $X$ are also preserved 
by the \textit{$N$-phone cord move} where $N$ is the rack rank of $X$:

\[\includegraphics{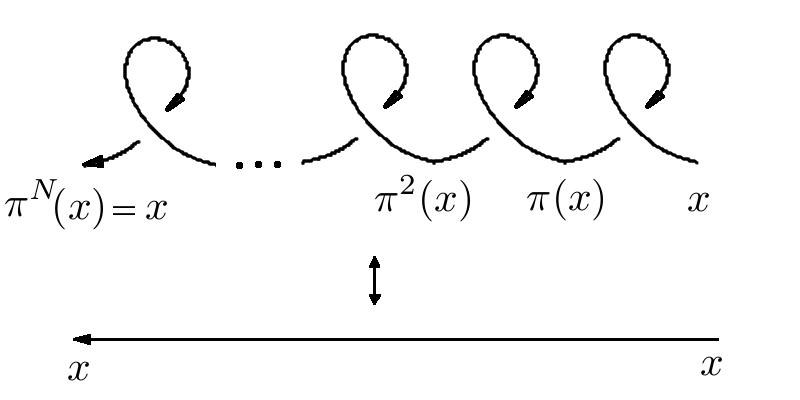}\]

In \cite{N} it is noted that if a rack $X$ has rack rank $N$, then the
basic counting invariants are periodic in $N$ on each component, and
thus the sum of these basic counting invariants over a complete set of
writhe vectors mod $N$ is an invariant of $L$ up to ambient isotopy.

\begin{definition}
\textup{The \textit{integral rack counting invariant} of $L$ with respect to 
$X$ is defined as:
\[\Phi^{\mathbb{Z}}_X(L) = \sum_{\mathbf{w}\in W} 
|\mathrm{Hom}(FR(L,\mathbf{w}),X)|\]
where $X$ is a finite rack, $L$ is a link with $c$ components, $N=N(X)$ and 
$W=(\mathbb{Z}_N)^c$. }
\end{definition}

\begin{example}
\textup{Let us find the set of rack labelings of the Hopf link $H$ by the
constant action rack $X_{(12)}=\{1,2\}$ with $\sigma=(12)$. The labeling 
rule here says that when crossing under a strand with any label, a 
$1$ switches to a $2$ and vice-versa. In particular, in this rack we always
have $x\ne x\tr y$ for all $x,y$. The rack rank is
$2$, so we need to find rack labelings over diagrams of $L$ with all 
writhe vectors in $(\mathbb{Z}_2)^2=\{(0,0),(0,1),(1,0),(1,1)\}$. The
diagrams of the Hopf link with writhe vectors
$(0,0)$, $(0,1)$ and $(1,0)$ have no valid rack labelings by $X$, since
each of these diagrams requires a relation of the form $x=x\tr y$:}
\[\begin{array}{ccc}
\includegraphics{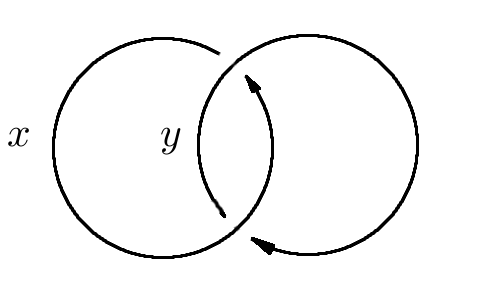} & \includegraphics{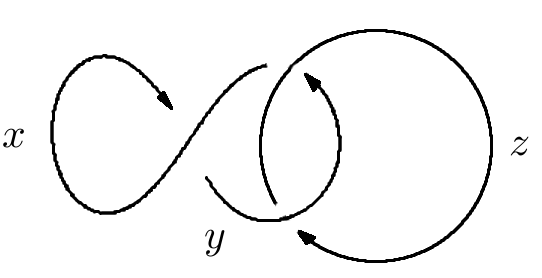} 
& \includegraphics{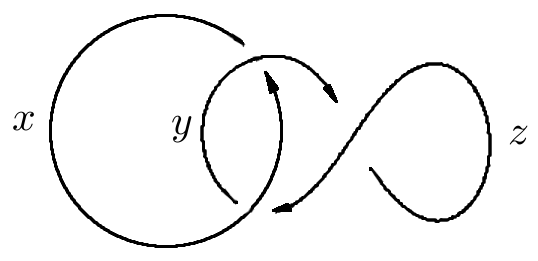} \\ 
x=x\tr y & z=z\tr y & x=x\tr y \\
y=y\tr x & x=y\tr z & y=z\tr x \\
         & y=x\tr x & z=y\tr z \\
\end{array}\]
\textup{while the $(1,1)$ diagram has four valid $X$-labelings:}
\[\scalebox{0.75}{\includegraphics{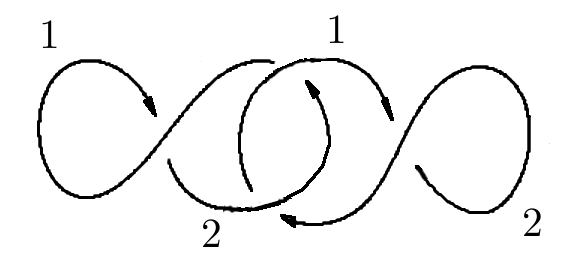} \ 
\includegraphics{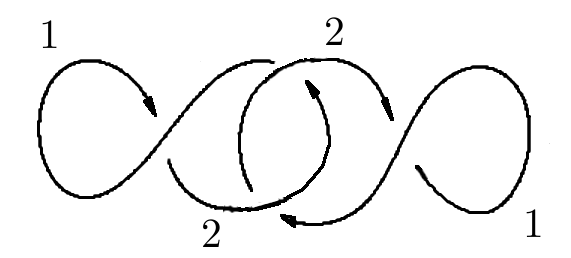} \ 
\includegraphics{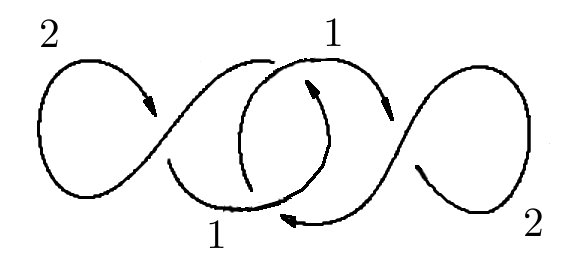} \ 
\includegraphics{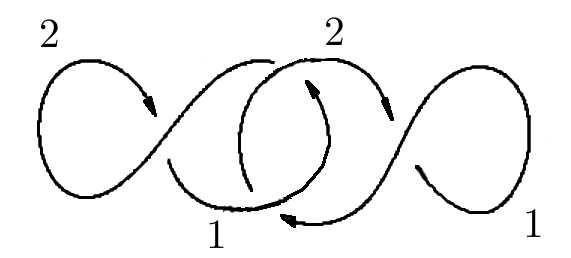}.} 
\]
\textup{Hence, the integral rack counting invariant of the Hopf link with 
respect to $X$ is $\Phi^{\mathbb{Z}}_X(H)=0+0+0+4=4$.}
\end{example}

\section{\large \textbf{The Rack Algebra and Rack Modules}}\label{ram}

In this section we recall and slightly reformulate the notions
of the rack algebra and rack modules from \cite{AG,CEGS,NJ}.

Let $X$ be a finite rack with rack rank $N$. We would like to define
an associative algebra $\mathbb{Z}[X]$ determined by $X$. One way to do 
this is to modify the $(t,s)$-rack structure described in section \ref{rb},
giving each pair $(x,y)\in X\times X$ its own invertible $t_{x,y}$ and
generic $s_{x,y}$. Let $\Omega[X]$ be the free $\mathbb{Z}$-algebra 
generated by these $s_{x,y}$ and $t_{x,y}^{\pm 1}$. The \textit{rack algebra}
of $X$, $\mathbb{Z}[X]$, will be a quotient of $\Omega[X]$ by a certain
ideal $I$.

To see the multiplicative structure of $\mathbb{Z}[X]=\Omega[X]/I$, it is
helpful to see a geometric interpretation. Let $L$ be a fixed knot or link
diagram with a fixed rack labeling by $X$. We would like to consider secondary
labelings of $L$ by elements of the rack algebra $\mathbb{Z}[X]$. In 
\cite{CEGS} these secondary labels are pictured as ``beads'' on the strands
of the knot or link. The basic bead-labeling rule is the usual $(t,s)$-rack
rule, with the exception that the beads at a crossing use the specific 
$t_{x,y}$ and $s_{x,y}$ associated to the rack labels on the right-hand 
underarc and overarc as illustrated.
\[\includegraphics{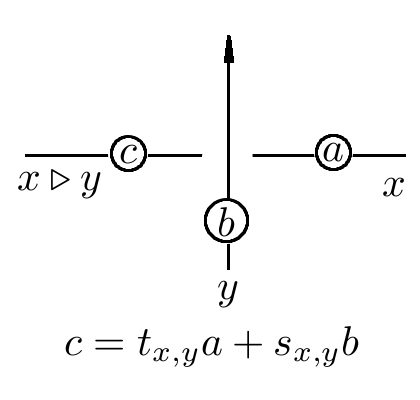}\]

\begin{remark}\textup{
Note that this rule does not make use of the orientation of the understrand. 
With rack label $x$ on the underarc on the right-hand side when the overarc 
is oriented upward and rack label $y$ on the overstrand, the rack algebra 
labels are $a$ on the right-hand underarc, $b$ on the overarc and  
$c=t_{x,y}a+s_{x,y}b$ on the left-hand overarc.}

\[\includegraphics{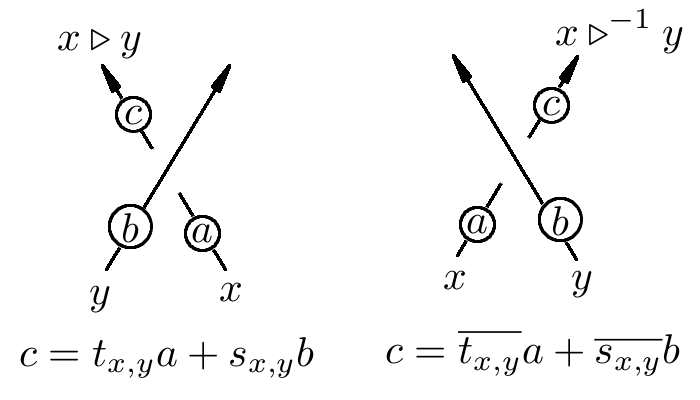}\]

\textup{We can also (as in \cite{CEGS}) formulate the rack algebra labeling 
rule in terms of the inbound arcs at a negative crossing; first we define 
elements $\bar{t_{x,y}}=t^{-1}_{x\tr^{-1}y,y}$ and 
$\bar{s_{x,y}}=-t^{-1}_{x\tr^{-1}y,y}s_{x\tr^{-1}y,y}$. Then we can express the 
label on the outbound underarc at a negative crossing with inbound
rack colors $x,y$ as
\[c=\bar{t_{x,y}}a+\bar{s_{x,y}}b
=t^{-1}_{x\tr^{-1}y,y}a-t^{-1}_{x\tr^{-1}y,y}s_{x\tr^{-1}y,y}b.\]
}\end{remark}

We would like to find the conditions required to make labelings of a 
rack-labeled diagram by rack algebra elements invariant under rack-labeled
blackboard framed isotopy. We will then define the rack algebra 
$\mathbb{Z}[X]$ by taking a quotient by the ideal $I$ generated by the 
elements we must kill to obtain invariance under $X$-colored blackboard 
framed Reidemeister moves. In \cite{AG,CEGS}
the special case where $X$ is a quandle, i.e. a rack of rack rank $N(X)=1$, 
is considered.  To extend this definition to the case of racks with finite 
rank $N\ge 1$, the Reidemeister moves we need are the standard type II 
and III moves together with the framed type I moves and the $N$-phone 
cord move. 

The Reidemeister II and blackboard framed Reidemeister I moves do not impose 
any conditions due to our choice of crossing rule:
\[\begin{array}{ccc}
\includegraphics{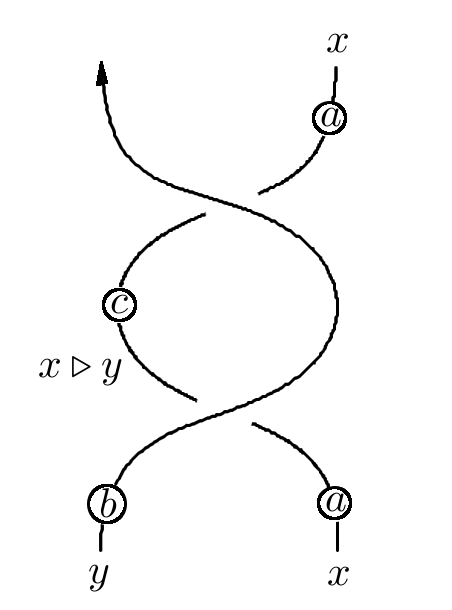} &
\includegraphics{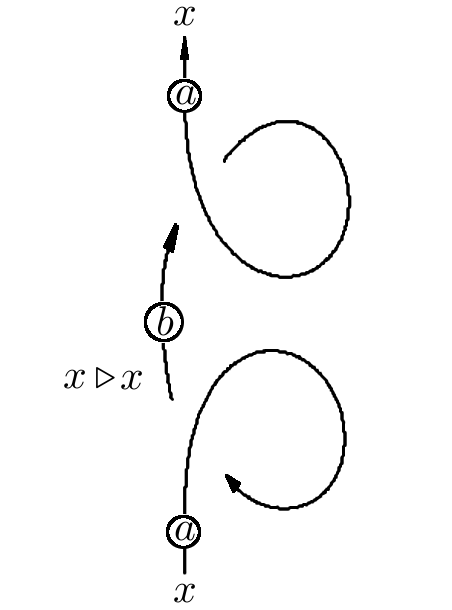} &
\includegraphics{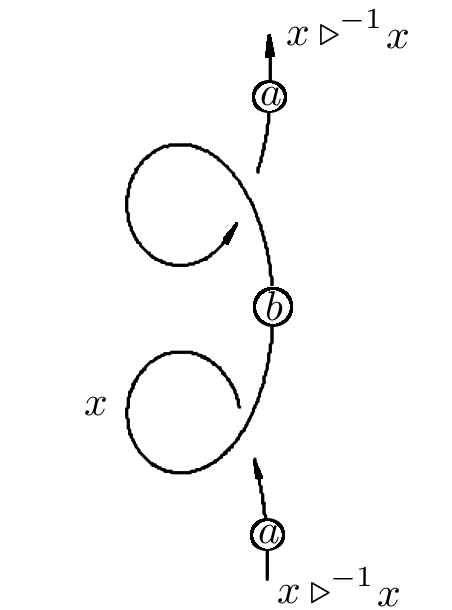}  \\
c=t_{x,y}a+s_{x,y}b &
b=(t_{x,x}+s_{x,x})a &
b=t_{x\tr^{-1}x,x}b+s_{x\tr^{-1}x,x}a \\
\end{array}
\]

The Reidemeister III move gives us three conditions which we recognize as
indexed versions of the $(t,s)$-rack conditions:
\[\includegraphics{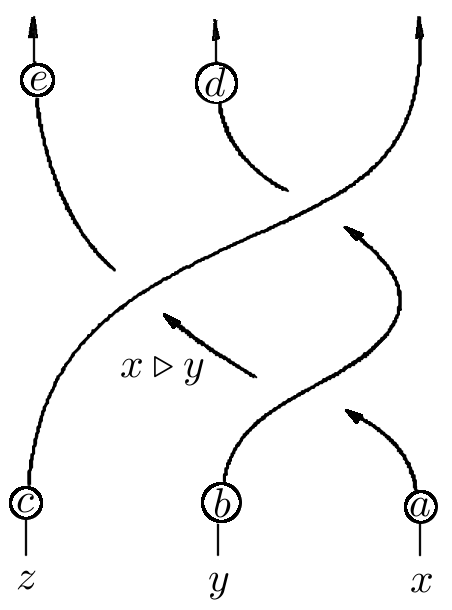}\quad\raisebox{1in}{$\sim$}\quad
\raisebox{-0.1in}{\scalebox{1.05}{\includegraphics{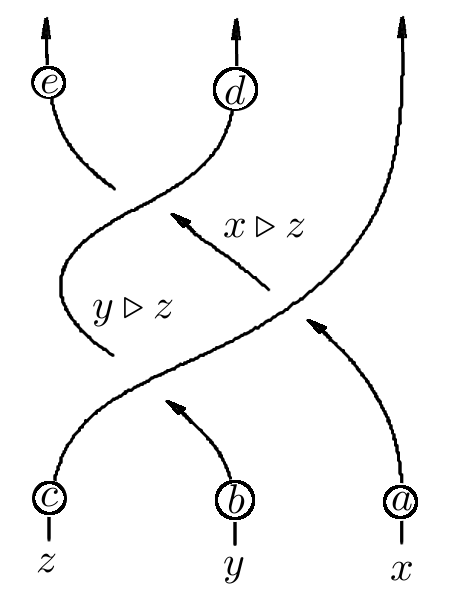}}}\]
Comparing coefficients on $e$, we must have
\[t_{x\tr y,z}t_{x,y}=t_{x\tr z, y\tr z}t_{x,z},
\quad t_{x\tr y,z}s_{x,y}=s_{x\tr z,y\tr z}t_{y,z}\]
\[ \mathrm{and} \quad s_{x\tr y,z}=s_{x\tr z,y\tr z}s_{y,z}+t_{x\tr z,y\tr z}s_{x,z}.\]
Finally, for defining counting invariants we will need $N$-phone cord
move compatibility. 
\[\includegraphics{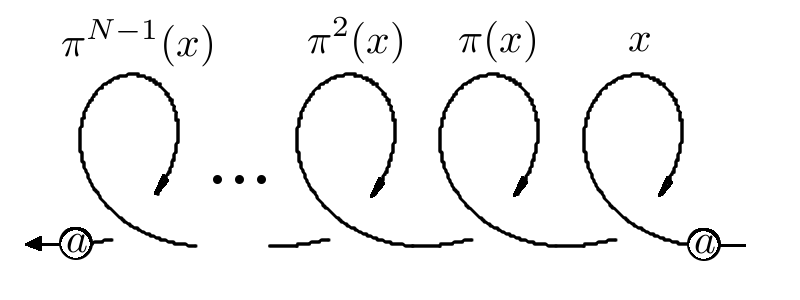}\]
Going though the $k$th kink multiplies our rack algebra element by a factor of
$t_{\pi^{k}(x),\pi^{k}(x)}+s_{\pi^{k}(x),\pi^{k}(x)}$, so $N$-phone cord invariance
requires that the product of these factors must equal 1:
\[\prod_{k=0}^{N-1}(t_{\pi^k(x),\pi^k(x)}+s_{\pi^k(x),\pi^k(x)}) =1.\]

Thus, we have:
\begin{definition}
\textup{Let $X$ be a rack with rack rank $N$ and let $\Omega[X]$ be the
free $\mathbb{Z}$-algebra generated by elements $t_{x,y}^{\pm 1}$ and $s_{x,y}$
where $x,y\in X$. Then the \textit{rack algebra} of $X$, denoted 
$\mathbb{Z}[X]$, is $\mathbb{Z}[X]=\Omega[X]/I$ where $I$ is the ideal 
generated by elements of the form}
\begin{list}{$\bullet$}{}
\item{$t_{x\tr y,z}t_{x,y}-t_{x\tr z, y\tr z}t_{x,z}$,}
\item{$t_{x\tr y,z}s_{x,y}-s_{x\tr z,y\tr z}t_{y,z}$,}
\item{$s_{x\tr y,z}-s_{x\tr z,y\tr z}s_{y,z}-t_{x\tr z,y\tr z}s_{x,z}$ \ \textup{and}}
\item{$\displaystyle{1-\prod_{k=0}^{N-1} (t_{\pi^k(x),\pi^k(x)}+s_{\pi^k(x),\pi^k(x)})}$}
\end{list}
\textup{A \textit{rack module over $X$} or a \textit{$\mathbb{Z}[X]$-module}
is a representation of $\mathbb{Z}[X]$, i.e., an abelian group $G$ with 
isomorphisms $t_{x,y}:G\to G$ and endomorphisms $s_{x,y}:G\to G$ such that 
each of the above maps is zero.}
\end{definition}

\begin{example}
\textup{Let $X$ be a rack and let $Y$ be a $(t,s)$-rack. Then $Y$ is
a $\mathbb{Z}[X]$-module via $t_{x,y}=t$ and $s_{x,y}=s$ for all $x,y$
if and only if $N(X)=N(Y)$.}
\end{example}

\begin{remark}\label{rem:ext}
\textup{Let $X$ be a rack of rack rank $N$ and let $G$ be an abelian group 
with isomorphisms $t_{x,y}:G\to G$ and endomorphisms $s_{x,y}:G\to G$. The 
rack module condition is equivalent to the condition that the cartesian 
product $X\times G$ is a rack of rack rank $N$ under the operation}
\[(x,g)\tr (y,h)=(x\tr y, t_{x,y}(g)+s_{x,y}(h)).\]
\textup{Such a rack is usually known as an \textit{extension rack} of $X$.}
\end{remark}

For the purpose of enhancing the rack counting invariant, we want
to find examples of finite rack modules. Probably the simplest way to 
achieve this is the following: let $X$ be a finite rack and let $R$ be a 
finite ring with identity, e.g. $\mathbb{Z}_n$ or, for a non-commutative 
example, the ring $M_{m}(\mathbb{Z}_n)$ of $m\times m$ matrices over 
$\mathbb{Z}_n$. We can give any $R$-module $V$ the structure of a 
$\mathbb{Z}[X]$-module by choosing elements $t_{x,y}\in R^{\ast}$ and 
$s_{x,y}\in R$ satisfying the rack module relations for all $x,y\in X$. 

We can express a $\mathbb{Z}[X]$-module structure on $R$ conveniently by 
giving a block matrix $M_{R}=[M_t|M_s]$ with block matrices $M_t$ and $M_s$ 
indexed by elements of $X=\{x_1,\dots, x_n\}$, i.e. the elements in row $i$ 
column $j$ of $M_t$ and $M_s$ respectively are $t_{x_i,y_j}$ and $s_{x_i,y_j}$.

\begin{example}\label{ex1}
\textup{Let $X$ be the constant action rack on the set $\{1,2\}$ with 
permutation $\sigma=(12)$, let $R=\mathbb{Z}_3$ and set  $s_{11}=s_{22}=1,$
$s_{12}=s_{21}=2$ and $t_{11}=t_{12}=t_{21}=t_{22}=1$. It is straightforward
to verify that the rack algebra relations are satisfied, and thus we have
a finite $\mathbb{Z}[X]$-module structure on $\mathbb{Z}_3$ with rack module 
matrix given by}
\[M_R=\left[\begin{array}{cc|cc}
1 & 1 & 1 & 2 \\
1 & 1 & 2 & 1
\end{array}\right].\]
\end{example}

\begin{example}
\textup{Let $L$ be a fixed oriented link diagram with arcs $a_1,\dots, a_n$ 
and let $f$ be a rack labeling of $L$ by $X$. We can regard 
$f\in\mathrm{Hom}(FR(L),X)$ as a rack homomorphism from the fundamental 
rack of $L$ into $X$ or as a link diagram with a fixed rack labeling. Such 
a labeling determines a rack module 
$\mathbb{Z}[f]=\mathbb{Z}[X][a_1,\dots,a_n]/C$ where $C$ is the submodule
determined by the crossing relations at the crossings of $L$. For instance,
the trefoil knot diagram below with the pictured rack coloring defines a 
rack module with the listed presentation matrix $M_{\mathbb{Z}[f]}$.}
\[\raisebox{-0.5in}{\includegraphics{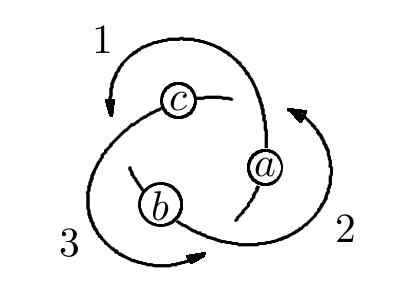}}
\quad
M_{\mathbb{Z}[f]}=\left[\begin{array}{rrr}
t_{13} & -1 & s_{13} \\
-1 & s_{32} & t_{32} \\
s_{21} & t_{21} & -1
\end{array}\right]
\]\end{example}

\begin{remark}
\textup{We note that in the case where $X$ is the trivial quandle with
one element, we have $t_{1,1}=t$, $s_{1,1}=s$ and $(t+s)=1$ implies $s=1-t$,
and $M_{\mathbb{Z}[f]}$ is a presentation matrix for the Alexander quandle
of $L$.}
\end{remark}

\section{\large \textbf{Rack Module Enhancements of the Counting Invariant}}
\label{rinv}

In this section we use finite rack modules to enhance the rack counting
invariant.

Let $X$ be a finite rack with rack rank $N$, $L$ an oriented link of 
$c$ ordered components, and $W=(\mathbb{Z}_N)^c$ the set of writhe
vectors of $L$ mod $N$. Recall from section \ref{rb} that the 
\textit{integral rack counting invariant} of $L$ with respect to $X$
is the sum
\[\Phi^{\mathbb{Z}}_{X}(L)=\sum_{\mathbf{w}\in W} 
|\mathrm{Hom}(FR(L,\mathbf{w}),X)|\]
counting the total number of rack labelings of $L$ by $X$ over a complete
period of writhes mod $N$.

Now let $R$ be a $\mathbb{Z}[X]$-module. For each rack labeling 
$f\in \mathrm{Hom}(FR(L,\mathbf{w}),X)$ of $L$ by $X$, we can ask how many 
associated $R$-labelings of the diagram there are. Such a labeling is 
a homomorphism of $\mathbb{Z}[X]$-modules $\phi:\mathbb{Z}[f]\to R$. 
We can find the set of all such labelings from the presentation matrix of 
$\mathbb{Z}[f]$ by substituting in the values of $t_{x,y}$ and $s_{x,y}$ 
in $R$ and finding the solution space of the resulting homogeneous system 
of linear equations in $R$. Denote this set by 
$\mathrm{Hom}(\mathbb{Z}[f],R)$.

By construction, the set $\mathrm{Hom}(\mathbb{Z}[f],R)$ of $R$-labelings of 
$f$ is invariant under $X$-colored blackboard-framed Reidemeister moves and
$N$-phone cord moves, and thus gives us an invariant signature of the 
$X$-labeling. The multiset of such signatures over the set of all $X$-labelings 
$\{\mathrm{Hom}(FR(L,\mathbf{w}),X)
\ :\ \mathbf{w}\in W\}$ gives us a natural enhancement of the rack counting 
invariant. 

\begin{definition}
\textup{Let $X$ be a finite rack with rack rank $N$, $L$ an oriented link of 
$c$ ordered components, $W=(\mathbb{Z}_N)^c$  and  $R$ a finite 
$\mathbb{Z}[X]$-module. The \textit{rack module enhanced multiset} of $L$ 
with respect to $X$ and $R$ is the multiset}
\[\Phi^{M}_{X,R}(L)=\{\mathrm{Hom}(\mathbb{Z}[f],R)\ 
:\ f\in \mathrm{Hom}(FR(L,\mathbf{w}),X), \mathbf{w}\in W\}\]
\textup{For ease of comparison, we can take the cardinality of each of the
sets of rack module homomorphisms to get a multiset of integers, the 
\textit{rack module enhanced counting multiset}}
\[\Phi^{M,\mathbb{Z}}_{X,R}(L)=\{|\mathrm{Hom}(\mathbb{Z}[f],R)| \ 
:\ f\in \mathrm{Hom}(FR(L,\mathbf{w}),X), \mathbf{w}\in W\}\]
\textup{The generating function of $\Phi^{M,\mathbb{Z}}_{X,R}(L)$ is a 
polynomial link invariant, the \textit{rack module enhanced counting invariant}}
\[\Phi_{X,R}(L)=\sum_{\mathbf{w}\in W} \left( 
\sum_{f\in \mathrm{Hom}(FR(L,\mathbf{w}),X)} u^{|\mathrm{Hom}(\mathbb{Z}[f],R)|}\right).
\]
\end{definition}

\begin{example}\label{ex2}
\textup{Let $X$ and $R$ be the rack and rack module in example \ref{ex1}. 
Every knot has exactly two labelings by $X$ in 
even writhe and zero labelings in odd writhe, so the integral rack counting
invariant for any knot is $\Phi_{X}^{\mathbb{Z}}(K)=2$. However, this smallest 
possible non-quandle rack with a choice of $\mathbb{Z}[X]$-module structure on
$\mathbb{Z}_3$ nevertheless detects the difference between the two smallest 
non-trivial knots, the trefoil $3_1$ and the figure-eight knot $4_1$.}

\textup{The $X$-colored trefoil knot below has the listed 
$\mathbb{Z}[f]$ presentation matrix.} 
\[\raisebox{-0.75in}{\includegraphics{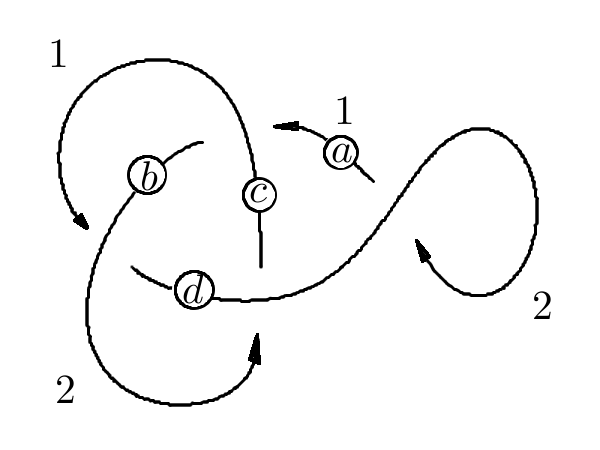}}  \quad
M_{\mathbb{Z}[f]}=
\left[\begin{array}{rrrr}
-1 & 0 & 0 & t_{2,2}+s_{2,2} \\
t_{1,1} & -1 & s_{1,1} & 0 \\
0 & s_{1,2} & t_{1,2} & -1 \\
0 & t_{2,2} & -1 & s_{2,2}\\
\end{array}\right]\]
\textup{Substituting in the values of $t_{x,y}$ and $s_{x,y}$ in $R$ and
thinking of the matrix as a coefficient matrix for a homogeneous system, 
we find there are $9$ total rack algebra labelings; note that the symmetry 
in the indices in $M_R$ implies that the other $X$-labeling yields the 
same number of $R$-labelings, for an invariant value of $\phi_{X,R}(3_1)=2u^9$}
\[\left[\begin{array}{rrrr}
2 & 0 & 0 & 2 \\
1 & 2 & 1 & 0 \\
0 & 2 & 1 & 2 \\
0 & 1 & 2 & 1\\
\end{array}\right] \rightarrow
\left[\begin{array}{rrrr}
1 & 0 & 0 & 1 \\
0 & 1 & 2 & 1 \\
0 & 0 & 0 & 0 \\
0 & 0 & 0 & 0
\end{array}\right]
\]

\textup{Comparing this with the $X$-colored figure eight knot $4_1$, we 
find that the space of $\mathbb{Z}[X]$-colorings has only three elements:}
\[\raisebox{-0.75in}{\includegraphics{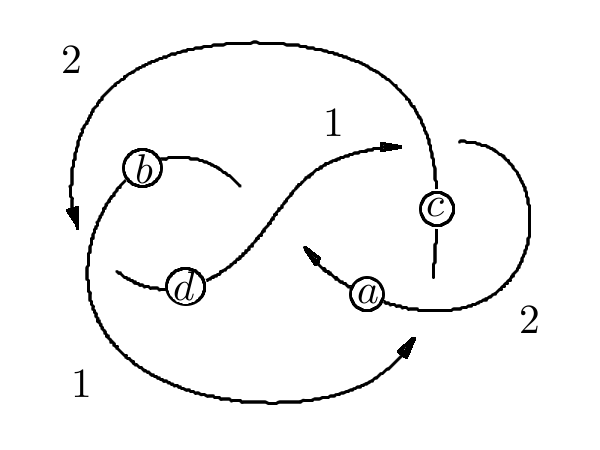}}  \quad
M_{\mathbb{Z}[f]}=
\left[\begin{array}{rrrr}
t_{2,2} & 0 & s_{2,2} & -1 \\
t_{2,1} & -1 & 0 & s_{2,1} \\
0 & s_{2,1} & t_{2,1} & -1 \\
s_{2,2} & -1 & t_{2,2} & 0\\
\end{array}\right]\]
\[\left[\begin{array}{rrrr}
1 & 0 & 1 & 2 \\
1 & 2 & 0 & 2 \\
0 & 2 & 1 & 2 \\
1 & 2 & 1 & 0 \\
\end{array}\right] \rightarrow
\left[\begin{array}{rrrr}
1 & 0 & 0 & 1 \\
0 & 1 & 0 & 2 \\
0 & 0 & 1 & 1 \\
0 & 0 & 0 & 0 \\
\end{array}\right].
\]
\textup{Thus, we have $\Phi_{X,R}(4_1)=2u^3\ne 2u^9=\Phi_{X,R}(3_1)$, and
the rack module enhanced invariant is strictly stronger than the rack 
counting invariant.}
\end{example}

\begin{remark}\textup{
In the special case when $R$ is a commutative ring
and $t_{x,y}$ and $s_{x,y}$ are elements of $R$, we do not need 
$R$ to be finite in order to have a well-defined computable enhanced invariant.
Replacing the the number of $R$-labelings with the dimension of the 
space of $R$-labelings yields a computable invariant we might call 
the \textit{dimension-enhanced rack module counting invariant}:
\[\Phi^{\mathrm{dim}}_{X,R}(L)=\sum_{\mathbf{w}\in W} \left( 
\sum_{f\in \mathrm{Hom}(FR(L,\mathbf{w}),X)} 
u^{\mathrm{dim}(\mathrm{Hom}(\mathbb{Z}[f],R))}\right)\]
Indeed, if $R$ is a PID, $\Phi^{\mathrm{dim}}_{X,R}$ carries the same 
information as $\Phi_{X,R}$, though if $R$ has torsion then $\Phi_{X,R}$ 
will generally contain more information.}
\end{remark}

\begin{remark}
\textup{Let $X$ be a rack and $R$ a finite $\mathbb{Z}[X]$-module. As
noted in remark \ref{rem:ext}, the cartesian product $X\times R$ is
a rack under}
\[(x,a)\tr(y,b)=(x\tr y, t_{x,y}a+s_{x,y}b).\]
\textup{The rack module enhanced counting invariant $\Phi_{X,R}$ is then
related to the integral rack counting invariant with respect to the 
extension rack $X\times R$ by} 
\[\Phi_{X\times R}^{\mathbb{Z}}(L)=\sum_{H\in\Phi^{M}_{X,R}(L)} |H|,\]
\textup{i.e., $\Phi_{X\times R}^{\mathbb{Z}}(L)$ is the sum of the exponents
in the terms of $\Phi_{X,R}.$ In particular, we can understand $\Phi_{X,R}$ 
as an enhancement of $\Phi_{X\times R}^{\mathbb{Z}}$ defined by collecting 
together the labelings of $L$ by the extension rack $X\times R$ which 
project to the same rack labeling of $L$ by $X$.}
\end{remark}

\begin{remark}
\textup{Note that all of the preceding extends to virtual knots and links
in the usual way, i.e., by ignoring any virtual crossings. See \cite{K} 
for more.}
\end{remark}

\section{\large \textbf{Computations and Examples}}\label{comp}

In this section we collect a few examples of the rack module enhanced 
invariant.

\begin{example}\label{ex3}
\textup{Let $X$ and $R$ be the rack and rack module used in examples 
\ref{ex1} and \ref{ex2}. The virtual knots labeled $3_1$ and $3_7$
in the Knot Atlas \cite{K} both have a 
Jones polynomial of 1. Since the rack $X$ has a rack rank of $2$ and each 
knot has writhe $1$, we must also account for colorings of each knot with 
the addition of a positive kink. In fact, the only valid colorings of both 
virtual knots $3_1$ and $3_7$ by the rack are the colorings of the knots 
with even writhe. One possible $X$-coloring of the virtual knot $3_1$ with 
its presentation matrix is:}
\[\raisebox{-0.75in}{\includegraphics{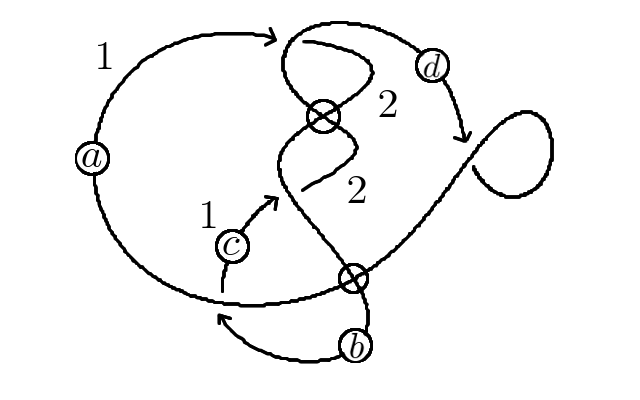}} 
\quad M_{\mathbb{Z}[f]}=\left[\begin{array}{cccc}
t^{-1}_{2,2} & -1 & 0 & -t^{-1}_{2,2}s_{2,2} \\
0 & s_{1,2} & t_{1,2} & -1 \\
-t^{-1}_{1,1}s_{1,1} & t^{-1}_{1,1} & -1 & 0 \\
s_{2,1}-1 & 0 & 0 & t_{2,1}
\end{array}\right].\]
\textup{Substituting the values from the rack module $R$ yields a matrix that 
may be reduced as follows:}
\[\left[\begin{array}{cccc}
1 & 2 & 0 & 2 \\
0 & 2 & 1 & 2 \\
2 & 1 & 2 & 0 \\
1 & 0 & 0 & 1
\end{array}
\right]\rightarrow
\left[\begin{array}{cccc}
1 & 1 & 0 & 0 \\
0 & 1 & 1 & 0 \\
0 & 0 & 1 & 1 \\
0 & 0 & 0 & 0
\end{array}
\right].\]
\textup{This row reduced matrix shows that there are three possible 
rack module colorings. In addition, the other possible coloring of this 
knot is symmetric to this coloring so the $3_1$ virtual knot has a rack 
module enhanced counting invariant of $2u^3$.}

\textup{In contrast, the $3_7$ virtual knot has presentation matrix}
\[M_{\mathbb{Z}[f]}=\left[\begin{array}{cccc}
t^{-1}_{2,1} & -1 & -t^{-1}_{2,1}s_{2,1} & 0 \\
-1 & -t^{-1}_{1,2}s_{1,2} & 0 & t^{-1}_{1,2} \\
s_{2,1} & t_{2,1} & -1 & 0 \\
0 & 0 & t_{1,2} & s_{1,2}-1
\end{array}\right]\]
\textup{derived from the coloring:}
\[\includegraphics{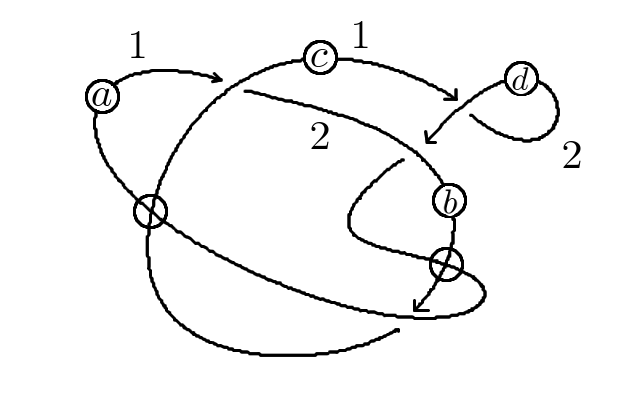}\]
\textup{The module values yield a presentation matrix that may be reduced 
as follows:}
\[\left[\begin{array}{cccc}
1 & 2 & 1 & 0 \\
2 & 1 & 0 & 1 \\
2 & 1 & 2 & 0 \\
0 & 0 & 1 & 1
\end{array}
\right]\rightarrow
\left[\begin{array}{cccc}
1 & 2 & 1 & 0 \\
0 & 0 & 1 & 1 \\
0 & 0 & 0 & 0 \\
0 & 0 & 0 & 0
\end{array}
\right].\]
\textup{This matrix indicates that the $3_7$ virtual knot with this 
coloring has $9$ rack module colorings. The $3_7$ knot also has one 
additional rack coloring, which is symmetric to the given coloring, so 
the rack module enhanced invariant for the $3_7$ virtual knot is $2u^9$. 
Both $3_1$ and $3_7$ have a Jones polynomial of 1, yet their rack module 
enhanced counting invariants are different, $2u^3 \ne 2u^9$, which shows 
that the the rack module enhanced invariant is not determined by the 
Jones polynomial.}
\end{example}

\begin{example}
\textup{The rack module enhanced counting invariant is also not determined 
by the Alexander polynomial. Knots $8_{18}$ and $9_{24}$ both have the same 
Alexander polynomial, $-t^3+5t^2-10t+13-10t^{-1}+5t^{-2}-t^{-3}$. Using the 
two element constant action rack and module used in examples \ref{ex1}, 
\ref{ex2}, and \ref{ex3}, the knot $8_{18}$ has two rack colorings, one of 
which is illustrated below.}
\[\scalebox{1}{\includegraphics{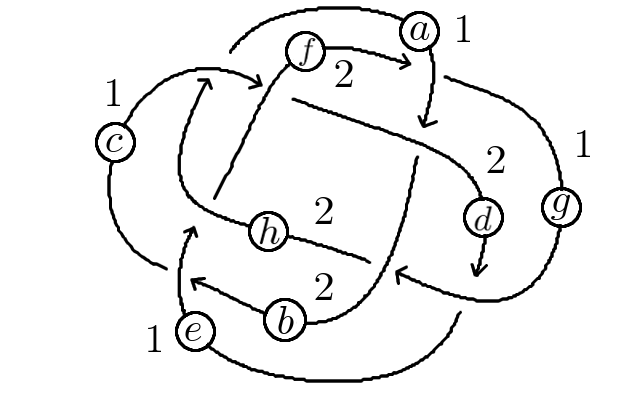}}\]
\textup{The presentation matrix of $8_{18}$ is:}
\[M_{\mathbb{Z}[f]}=\left[\begin{array}{cccccccc}
s_{2,1} & 0 & 0 & 0 & 0 & t_{2,1} & -1 & 0 \\
t^{-1}_{2,2} & -1 & 0 & -t^{-1}_{2,2}s_{2,2} & 0 & 0 & 0 & 0 \\
0 & -t^{-1}_{2,2}s_{2,2} & 0 & 0 & 0 & 0 & t^{-1}_{2,2} & -1 \\
0 & t_{2,1} & -1 & 0 & s_{2,1} & 0 & 0 & 0 \\
-1 & 0 & s_{2,1} & 0 & 0 & 0 & 0 & t_{2,1} \\
0 & 0 & t^{-1}_{2,2} & -1 & 0 & -t^{-1}_{2,2}s_{2,2} & 0 & 0 \\
0 & 0 & 0 & t_{2,1} & -1 & 0 & s_{2,1} & 0 \\
0 & 0 & 0 & 0 & t^{-1}_{2,2} & -1 & 0 & -t^{-1}_{2,2}s_{2,2}
\end{array}\right]\]
\textup{which can be reduced as shown after substituting in values from 
the module:}
\[\left[\begin{array}{cccccccc}
2 & 0 & 0 & 0 & 0 & 1 & 2 & 0 \\
1 & 2 & 0 & 2 & 0 & 0 & 0 & 0 \\
0 & 2 & 0 & 0 & 0 & 0 & 1 & 2 \\
0 & 1 & 2 & 0 & 2 & 0 & 0 & 0 \\
2 & 0 & 2 & 0 & 0 & 0 & 0 & 1 \\
0 & 0 & 1 & 2 & 0 & 2 & 0 & 0 \\
0 & 0 & 0 & 1 & 2 & 0 & 2 & 0 \\
0 & 0 & 0 & 0 & 1 & 2 & 0 & 2
\end{array}\right]\rightarrow \left[\begin{array}{cccccccc}
1 & 2 & 0 & 2 & 0 & 0 & 0 & 0 \\
0 & 1 & 2 & 0 & 2 & 0 & 0 & 0 \\
0 & 0 & 1 & 2 & 2 & 0 & 0 & 0 \\
0 & 0 & 0 & 1 & 1 & 1 & 2 & 1 \\
0 & 0 & 0 & 0 & 1 & 2 & 0 & 2 \\
0 & 0 & 0 & 0 & 0 & 0 & 0 & 0 \\
0 & 0 & 0 & 0 & 0 & 0 & 0 & 0 \\
0 & 0 & 0 & 0 & 0 & 0 & 0 & 0
\end{array}\right].\]
\textup{This matrix indicates that there are twenty-seven rack module 
colorings of the knot $8_{18}$, yielding an invariant value of $2u^{27}$ 
since there is one additional rack coloring due to the symmetry of the 
rack and rack module used.}

\textup{The knot $9_{24}$ requires a writhe adjustment by adding one 
positive kink in order to have a valid coloring by the rack, 
which yields the knot with presentation matrix shown below.}
\[\scalebox{1}{\includegraphics{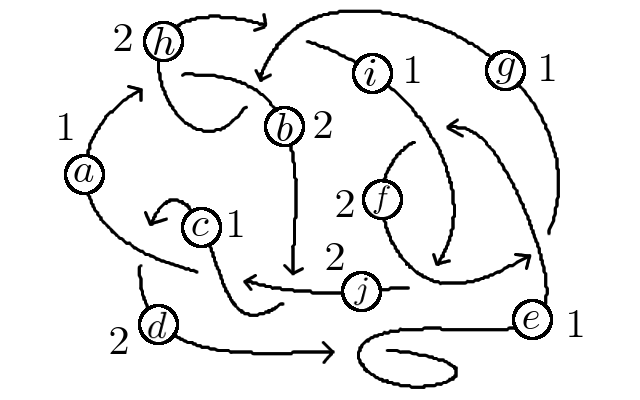}}\]
\[M_{\mathbb{Z}[f]}=\left[\begin{array}{cccccccccc}
t^{-1}_{2,2} & -1 & 0 & 0 & 0 & 0 & 0 & -t^{-1}_{2,2}s_{2,2} & 0 & 0 \\
0 & -t^{-1}_{2,2}s_{2,2} & 0 & 0 & 0 & 0 & t^{-1}_{2,2} & -1 & 0 & 0 \\
0 & t_{2,2} & -1 & 0 & 0 & 0 & 0 & 0 & 0 & s_{2,2} \\
-1 & 0 & s_{2,1} & 0 & 0 & 0 & 0 & 0 & 0 & s_{2,1}\\
s_{1,1} & 0 & t_{1,1} & -1 & 0 & 0 & 0 & 0 & 0 & 0\\
0 & 0 & 0 & t^{-1}_{1,1} & -t^{-1}_{1,1}s_{1,1}-1 & 0 & 0 & 0 & 0 & 0\\
0 & 0 & 0 & 0 & -t^{-1}_{1,1}s_{1,1} & t^{-1}_{1,1} & -1 & 0 & 0 & 0\\
0 & 0 & 0 & 0 & t^{-1}_{2,1} & -1 & 0 & 0 & -t^{-1}_{2,1}s_{2,1} & 0\\
0 & 0 & 0 & 0 & 0 & -t^{-1}_{2,2}s_{2,2} & 0 & 0 & t^{-1}_{2,2} & -1\\
0 & 0 & 0 & 0 & 0 & 0 & s_{2,1} & t_{2,1} & -1 & 0
\end{array}\right]\]
\textup{After substituting in the values from the rack module, the 
matrix can be reduced to the form:}
\[\left[\begin{array}{cccccccccc}
1 & 2 & 0 & 0 & 0 & 0 & 0 & 2 & 0 & 0 \\
0 & 2 & 0 & 0 & 0 & 0 & 1 & 2 & 0 & 0 \\
0 & 1 & 2 & 0 & 0 & 0 & 0 & 0 & 0 & 1 \\
2 & 0 & 2 & 0 & 0 & 0 & 0 & 0 & 0 & 1\\
1 & 0 & 1 & 2 & 0 & 0 & 0 & 0 & 0 & 0\\
0 & 0 & 0 & 1 & 1 & 0 & 0 & 0 & 0 & 0\\
0 & 0 & 0 & 0 & 2 & 1 & 2 & 0 & 0 & 0\\
0 & 0 & 0 & 0 & 1 & 2 & 0 & 0 & 1 & 0\\
0 & 0 & 0 & 0 & 0 & 2 & 0 & 0 & 1 & 2\\
0 & 0 & 0 & 0 & 0 & 0 & 2 & 1 & 2 & 0
\end{array}\right]\rightarrow \left[\begin{array}{cccccccccc}
1 & 2 & 0 & 0 & 0 & 0 & 0 & 2 & 0 & 0 \\
0 & 1 & 2 & 0 & 0 & 0 & 0 & 0 & 0 & 1 \\
0 & 0 & 1 & 0 & 0 & 0 & 2 & 1 & 0 & 2 \\
0 & 0 & 0 & 1 & 1 & 0 & 0 & 0 & 0 & 0\\
0 & 0 & 0 & 0 & 1 & 0 & 0 & 0 & 0 & 1\\
0 & 0 & 0 & 0 & 0 & 1 & 0 & 0 & 2 & 1\\
0 & 0 & 0 & 0 & 0 & 0 & 1 & 1 & 0 & 0\\
0 & 0 & 0 & 0 & 0 & 0 & 0 & 1 & 1 & 0\\
0 & 0 & 0 & 0 & 0 & 0 & 0 & 0 & 0 & 0\\
0 & 0 & 0 & 0 & 0 & 0 & 0 & 0 & 0 & 0
\end{array}\right].\]
\textup{The reduced matrix shows that there are nine colorings of this 
knot by the rack module. The other coloring of the knot by the rack yields 
an additional nine colorings so the rack module enhanced invariant for 
the knot $9_{24}$ is $2u^9$. Since $2u^{27}\ne 2u^9$, the rack module 
enhanced counting invariant is not determined by the Alexander polynomial.} 
\end{example}

Our next examples give a sense of the effectiveness of some examples of 
$\Phi_{X,R}$ using racks $X$ and rack modules $R=\mathbb{Z}_n$ chosen at 
random from the results of our computer searches.

\begin{example} \label{ex:a1}
\textup{Let $X$ be the quandle with quandle matrix given by}
\[M_X=
\left[\begin{array}{rrrr}
1 & 3 & 2 \\
3 & 2 & 1 \\
2 & 1 & 3 \\
\end{array} \right] \]
\textup{Via computations in python, we selected a  $\mathbb{Z}[X]$-module 
structure on $R=\mathbb{Z}_5$ given by the rack module matrix}
\[M_R=\left[\begin{array}{ccc|ccc}
4 & 2 & 3 & 2 & 1 & 1 \\
1 & 4 & 2 & 2 & 2 & 3 \\
1 & 3 & 4 & 3 & 2 & 2
\end{array}\right] \]
\textup{and computed the knot invariant for all prime knots with up to 8 
crossings, all prime links with up to 7 crossings as listed in the 
Knot Atlas \cite{KA}, the unknot $U$, the 
square knot $SK$ and granny knot $GK$, and unlinks of two and three components
$U_2$ and $U_3$. The results are collected in the table below.}

\begin{center}
\begin{tabular}{r|l}
$\Phi_{X,M}$ & $L$ \\ \hline
$3u^5$ & $U, 5_2,6_2, 6_3, 7_1, 7_2, 7_3, 7_5, 7_6,8_1, 8_2, 8_3, 8_4, 8_6, 8_7, 8_{12}, 8_{13}, 8_{14}, 8_{17}$ \\ 
 & $L2a1, L4a1, L5a1, L6a4, L6n1, 7La3, L7a4, L7a6, L7n1, L7n2$ \\
$3u^5+6u^{25}$ & $3_1, 6_1, 7_7, 8_5, 8_{10}, 8_{11}, 8_{15}, 8_{19}, 8_{20}, 
L6a1, L6a3, L6a5, L7a1, L7a5$ \\
$3u^5+12u^{25}+12u^{125}$ & $SK, GK$ \\
$3u^{25}$ & $4_1, 5_1, 8_8, 8_9, 8_{16}, L6a2, L7a2, L7a7$ \\
$9u^{25}$ & $7_4, U_2$ \\
$27u^{25}$ & $8_{18}, 8_{21}, 9_2, 9_{24}, U_3$
\end{tabular}
\end{center}
\end{example}

\begin{example} \textup{Let $X$ be the rack with rack matrix given by}
\[
M_X=\left[\begin{array}{rrrr}
2 & 2 & 2 \\
1 & 1 & 1 \\
3 & 3 & 3 \\
\end{array} \right] \]
\textup{Via computations in python, we selected a 
$\mathbb{Z}[X]$-module structure on $R=\mathbb{Z}_{5}$ given by}
\[M_{R}=\left[\begin{array}{ccc|ccc}
1 & 1 & 1 & 0 & 0 & 4 \\
1 & 1 & 1 & 0 & 0 & 4 \\
1 & 1 & 1 & 0 & 0 & 0
\end{array}\right] \]
\textup{and computed the knot invariant for the same set of knots and link
as in example \ref{ex:a1}. For all of the knots in our list, the 
invariant value is $\Phi_{X,M}=4u^5$. Our results for links are collected 
in the table below.}

\begin{center}
\begin{tabular}{r|l}
$\Phi_{X,M}$ & $L$ \\ \hline
$8u^5+8u^{25}$ & $L2a1,L4a1,L6a1,L6a2,L6a3,L7a2,L7a5,L7a6,L7n1$ \\
$16u^{25}$ & $U_2,L5a1,L7a1,L7a3,L7a4,L7n2$ \\
$48u^{25}+ 16u^{125}$ & $L6a4,L6a5,L6n1,L7a7$ \\
$64u^{125}$ & $U_3$ \\
\end{tabular}
\end{center}
\end{example}

Our \texttt{python} and \texttt{C} code for computing rack modules
and their invariants is available for download at \texttt{www.esotericka.org}.

\section{\large \textbf{Open Questions}}\label{q}

In this section we collect a few questions for future investigation.

\begin{itemize}
\item{If $X$ is a $(t,s)$-rack with rack rank $N=1$, i.e an Alexander quandle, 
then $X$ is a rack module over the trivial quandle with one element $X'$ 
with $t_{1,1}=t$ and $s_{1,1}=1-t$. In this scenario, the integral rack 
counting invariant $\Phi_{X}^{\mathbb{Z}}(L)$ is related to the rack module 
enhanced counting invariant $\Phi_{X',X}$  by
\[\Phi_{X',X}(L)=u^{\Phi_{X}^{\mathbb{Z}}(L)}.\]
How does this relationship generalize when $X'$ is a larger $(t,s)$-rack?}
\item{Other enhancements of the rack counting invariant using rack modules
are possible; we have only scratched the surface with the most basic
enhancement. What other enhancements of $\Phi^{\mathbb{Z}}_{X}$ using rack 
modules are possible?}
\item{What enhancements of $\Phi_{X,M}$ are possible?}
\item{What approach might we use to categorify these invariants?}
\item{How is $\Phi_{X,M}$ related to other knot and link invariants?}
\item{How do rack modules extend to the case of virtual racks, i.e.
racks with a nontrivial action at virtual crossings?}
\end{itemize}

\end{document}